\magnification1200
\baselineskip20pt

\def\ep{\varepsilon} 
\def\ch{\raise.4ex\hbox{$\chi$}}
\def\va{\raise.4ex\hbox{$\varphi$}}
\def\centreline{\centerline} \def\equalign{\eqalign}
  \def\lequalignno{\leqalignno}
\def\nt{\noindent} 
 \def\ga{\raise.4ex\hbox{$\gamma$}}

\def\Vv{{|||}}   \def\pf{\noindent{\bf Proof.\ \ }} 
    
\font\titlefont=cmbx10 scaled\magstep 3
\font\namefont=cmbx10 scaled\magstep 1
\def\ref{\noindent\hang}
\def\SA{\lower.9ex\hbox{$\buildrel \subset \over \rightarrow$}}
\def\BTD{\raise.4ex\hbox{$\bigtriangledown$}}

\def\MATHOP(#1){\mathop{#1}}

\def\AST{\lower.8ex\hbox{*}}

\def\h{{1\over 2}}
 3
 1
\def\Rogerslash{\backslash\kern-2pt\backslash}
\def\RS{/\kern-2pt /}

\def\GH{G{\mathop{\times
\atop{\raise.4ex\hbox{\sevenrm{\ \ s.d.} }}}}H}




\def\vz{\ |\kern-3pt-\kern-5pt-\ }
\def\vr{\ |\kern-2pt-\kern-5pt\rightarrow\ }

\def\ontoarrow{\geq\kern-3pt-\kern-5pt-\kern\gg}


\def\bookformat{ 
\hsize=157true mm\vsize=240true mm\parindent=25true mm
\baselineskip=18true pt \voffset=36true pt \nopagenumbers
\def\rhead{\hfil} \def\lhead{\hfil}
 \headline={\sevenrm\it\ifodd\pageno\rhead\else\lhead\fi}
\def\newsection##1{\vskip20true pt{\bf##1}
\def\rhead{\centreline{##1}\rlap{\mit\folio}}
\def\lhead{\llap{\mit\folio}\centreline{##1}}}
\def\newchapter##1{{\bf##1}\vskip53true mm}}

\def\newleftmargin#1mm#2mm{
\def\funnymargin{\ifodd\pageno#1mm\else#2mm\fi}
\hoffset=\funnymargin}
 
\def\onto{\ >\kern-5pt-\kern-2pt\kern-4pt\gg\ }

\input amssym.def
\input amssym

\centreline{\titlefont{The Compact Approximation
Property}}  

\medskip
\centreline{\titlefont{does not imply}} 

\medskip
\centreline{\titlefont{the Approximation Property}}

\vskip.5in
\centreline{\namefont{by}}

\bigskip
\centreline{\namefont{George Willis}}

\bigskip
\centreline{\bf{Mathematics Research Section}}

\centreline{\bf{School of Mathematical Sciences}}

\centreline{\bf{The Australian National University}}

\centreline{\bf{GPO Box 4,\qquad Canberra, ACT}}

\centreline{\bf{Australia}}

\vskip1in

\nt{\bf Abstract}:  It is shown how to construct, given a
Banach space which does not have the approximation
property, another Banach space which does not have the
approximation property but which does have the compact
approximation property.

\vskip1.5in
\nt{\it 1991 Mathematics Subject Classification: primary
46B28; secondary 46B10}

\vfill\eject
A Banach space, $X$, is said to have the
{\underbar{approximation}} {\underbar{property}} if for
every compact set $K \subseteq X$ and every $\ep > 0$
there is a finite rank operator $T$ on $X$ such that
$\Vert Tx - x\Vert < \ep$ for every $x$ in $K$.  The
approximation property is weaker than the notion of a
Schauder basis.  Classical Banach spaces have Schauder
bases but it was shown by Enflo [E] that there are spaces
which do not have the approximation property.  A shorter
construction of such spaces was given by Davie [D1],
[D2], see also [L \& T1], Theorem 2.d.3.

A Banach space, $X$, is said to have the
{\underbar{compact}} {\underbar{approximation}}
{\underbar{property}} if for every compact set $K
\subseteq X$ and every $\ep > 0$ there is a compact
operator $T$ on $X$ such that $\Vert Tx - x\Vert < \ep$
for every $x$ in $K$.  A finite rank operator is compact
and so the compact approximation property is formally
weaker than the approximation property.  Known examples
leave open the possibility that these two properties are
equivalent.  Indeed, many of the spaces which do not have
the approximation property are known to also not have the
compact approximation property, see [L \& T1] p.94 and
[L \& T2] Theorem 1.g.4.  However, the following
construction produces spaces which have the compact
approximation property while not having the approximation
property, thus showing that the two properties  are not
equivalent.

The construction is based on an argument due to 
Grothendieck [G] which shows that $X$ has the
approximation property if and only if for every Banach
space $Y$, every compact operator $T : Y \rightarrow X$
and every $\ep > 0$ there is a finite rank $F : Y
\rightarrow X$with $\Vert T - F\Vert < \ep$.  We shall
follow the exposition given in [L \& T1], Theorem 1.e.4.

Let $X$ be a Banach space which does not have the
approximation property.  Then there is a compact set $K
\subseteq X$ such that the identity operator cannot be
approximated on $K$ by finite rank operators.  By a
theorem of Grothendieck, see [L \& T1], Proposition 1.e.2,
it may be supposed that $K = {\overline{\rm conv}} \{
x_n\}^\infty_{n=1}$ where $\Vert x_n\Vert \leq 1$ for all
$n$ and $\Vert x_n\Vert$ decreases to zero.

For each $t$ between $0$ and $1$, put $U_t =
{\overline{\rm conv}} \{ \pm x_n/\Vert
x_n\Vert^t\}^\infty_{n=1}$.  Then $U_t$ is a compact,
convex, symmetric subset of $X$.  Let $Y_t$ be the linear
span of $U_t$ and define a norm on $Y_t$ by $\Vv x \Vv_t
= \inf \{\vert \lambda\vert : \lambda^{-1}\,x \in U_t\}$,
$x \in Y_t$.  It may be checked that
$(Y_t,\,\Vv\cdot\Vv_t)$ is a Banach space with unit ball
$U_t$. If $s < t$, then $U_s \subseteq U_t$ and so $Y_s
\subseteq Y_t$ and all spaces are contained in $X$. 
Denote the inclusion map of $Y_t$ into $X$ by $L_t$.  Then
$L_t$ is compact and has norm at most one.  It is shown in
[L \& T1], Theorem 1.e.4, that the operator $L_{1\over 2}
: Y_{1\over 2} \rightarrow X$ cannot be approximated by
finite rank operators.

The space to be constructed will be a space of functions
on $(0,\,1)$ with values in $X$.  For each $(s,\,t)
\subseteq (0,\,1)$ and $y$ in $Y_s$, $y\,{{\cal
X}_{(s,t)}}$ is such a function, where $y\,{{\cal
X}_{(s,t)}}$ denotes the map
$$r \mapsto\cases{y\ ,\quad &if $r \in (s,\,t)$\cr
0\ ,&otherwise\ \ \ .\cr}$$
Let $\cal Z$ be the linear span of $\{ y\,{{\cal
X}_{(s,t)}} : 0 < s < t < 1\ ;\ y \in Y_s\}$.  Note
that if $f$ belongs to $\cal Z$, then $f(r)$ is in $Y_r$
for all $r$.  Hence we may define a norm on $\cal Z$ by
$$\Vert f \Vert = \int^1_0 \Vv f(r)\Vv_r\,dr\ \ ,\ \ (f
\in {\cal Z})\ \ \ .$$
Now let $Z$ be the completion of $\cal Z$ with
respect to this norm.
  
\bigskip
\proclaim
Proposition 1.  $Z$ does not have the approximation
property.

\bigskip
\pf Define a map $R : Y_{1\over 2} \rightarrow Z$ by
$$R(y) = 2y\,{{\cal X}_{({1\over 2},1)}}\ \ \ ,\ \ \ (y
\in Y_{1\over 2})\ \ \ .\leqno(1)$$
Then
$$\lequalignno{\Vert R\,x_n\Vert &= 2 \int^1_{1\over
2}\,\Vv x_n\Vv_r\,dr\cr
{}&\leq 2 \int^1_{1\over 2} \Vert x_n\Vert^r\,dr\ \ \ ,\
\ \ {\rm because}\ \Vv x_n\Vv_r \leq \Vert x_n\Vert^r\ \
\ ,\cr
{}&< 2\,\Vert x_n\Vert^{1\over 2}/\big\vert \,{\rm
ln}\,\Vert x_n\Vert \,\big\vert\ \ \ .&(2)\cr}$$
Since $\Vert x_n\Vert \rightarrow 0$ as $n \rightarrow
\infty$, it follows that $\Vert R(x_n/\Vert
x_n\Vert^{1\over 2})\Vert \rightarrow 0$ as $n
\rightarrow \infty$, whence $R\,U_{1\over 2}$ is a totally
bounded subset of $Z$.  Since $U_{1\over 2}$ is the
unit ball in $Y_{1\over 2}$ it follows that $R$ is a
compact operator.

Now define a map $J : Z \rightarrow X$ by
$$J(f) = \int^1_0 f(r)\,dr\ \ \ ,\ \ \ (f \in Z)\ \ \ ,
\leqno(3)$$
where the integral may be defined in the obvious way if 
$f$ is in $\cal Z$ and this definition extends to all of
$Z$ by continuity.  Then $JR = L_{1\over 2}$ and $L_{1\over
2}$ cannot be approximated by finite rank operators.  It
follows that $R$ cannot be approximated by finite ranks. 
Therefore, by [L \& T1] Theorem 1.e.4, $Z$ does not have
the approximation property.{$\spadesuit$}

\bigskip
\proclaim
Proposition 2.  $Z$ does have the compact approximation
property.

\bigskip
\pf For each $r$ in $(0,\,1)$ define the operator which
shifts by $r$, $S_r$, on $Z$ as follows:  first, for
$(s,\,t) \subseteq (0,\,1)$ and $y$ in $T_s$, define 
$S_r(y\,{{\cal X}_{(s,t)}}) = y\,{{\cal X}_{(s+r,t+r)}}$;
next extend $S_r$ to $\cal Z$ by linearity; and then,
since $S_r$ is clearly a contraction mapping, extend to
$Z$ by continuity.  As already mentioned, $\Vert S_r\Vert
\leq 1$ for each $r$.  Furthermore, it may easily be
checked that the function $r \mapsto S_r\,f$ is norm
continuous, and $\Vert S_r\,f - f\Vert \rightarrow 0$ as
$r \rightarrow 0$ for each $f$ in $Z$.  However, $S_r$ is
not a compact operator.

To obtain compact operators on $Z$ which approximate the
identity define, for each $n$, $T_n : Z \rightarrow Z$, by
$$T_n\,f = n\,\int^{1/n}_0\,S_r\,f\,dr\ \ \ ,\ \ \ (f
\in Z)\ \ \ .\leqno(4)$$
The integral exists because $r \mapsto S_r\,f$ is norm
continuous and $\Vert T_n\Vert \leq 1$.  Since the
operators $S_r$ approximate the identity as $r$
approaches zero, it follows that $\Vert T_n\,f - f\Vert
\rightarrow 0$ as $n \rightarrow \infty$ for each $f$ in
$Z$.  Hence, since $\Vert T_n\Vert \leq 1$ for every
$n$, $T_n$ approximates the identity operator on a given
compact set when $n$ is sufficiently large.

We shall see that $T_n$ is compact for each $n$ by
showing that $T_n\,Z^{(1)}$ is totally bounded, where
$Z^{(1)}$ denotes the unit ball in $Z$.  First note that
the functions of the form
$$\sum^p_{i=1}\,\lambda_i (t_i - s_i)^{-1}\,y_i\,{\cal
X}_{(s_i,t_i)}$$
(where: $s_1 < t_1 < s_2 < t_2 < \ldots < s_p < t_p$;
$y_i \in Y_{s_i}$, $\Vv y_i\Vv_{s_i} \leq 1$; and
$\sum^p_{i=1}\,\vert \lambda _i\vert = 1$) are dense in
$Z^{(1)}$.  Hence, it will suffice to show that $T_n\{(t
- s)^{-1}\,y\,{\cal X}_{(s,t)} : s < t\ ;\ y \in U_s\}$
is totally bounded.  Next, since the unit ball in $Y_s$
is ${\overline{\rm conv}}\{\pm x_m/\Vert
x_m\Vert^s\}^\infty_{m=1}$, it will suffice to show that
$$T_n\,\{ (t - s)^{-1}\,\Vert x_m\Vert^{-s}\,x_m\,{\cal
X}_{(s,t)} : s < t\ ;\ m = 1, 2, 3, \ldots\}$$
is totally bounded.

For each $s$ and $t$ with $s < t$ we have that $x_m$
belongs to $Y_s$ for each $m$ and
$$T_n (x_m\,{\cal X}_{(s,t)}) = x_m\,h\ \ \ ,$$
where $h = n\,\int^{1/n}_0\,{\cal X}_{(s+r,t+r)}\,dr$.

For functions, $f$ and $g$, in $L^1(0,\,1)$, let $f * g$
denote the usual convolution product of $f$ and $g$
restricted to $(0,\,1)$.  Then $h = {\cal X}_{(s,t)} *
(n\,{\cal X}_{(0,1/n)})$.  It follows that for each $f$
in $L^1(0,\,1)$, $T_n (x_m\,f) = x_m (f * (n\,{\cal
X}_{(0,1/n)}))$.  Now it is well known, and may easily
be checked, that the map $f \mapsto f * (n\,{\cal
X}_{(0,1/n)})$ is a compact operator on $L^1(0,\,1)$. 
Hence for each $m$ the set $T_n\{(t - s)^{-1}\,\Vert
x_m\Vert ^{-s}\,x_m\,{\cal X}_{(s,t)} : s < t\}$ is totally
bounded. Furthermore,
$$\lequalignno{\Vert T_n (x_m\,{\cal X}_{(s,t)})\Vert &= n
\int^1_0 \Vv x_m\Vv_r \big\vert {\cal X}_{(s,t)} * {\cal
X}_{(0,1/n)} (r) \big\vert\,dr&(5)\cr
{}& \leq n\,\int^{t+1/n}_s \Vert x_m\Vert^r (t - s)\,dr\cr
{}&< n (t - s)\,\Vert x_m\Vert^s/\big\vert {\rm ln}\,\Vert
x_m\Vert\,\big\vert\ \ \ .\cr}$$
Hence, for each $\ep > 0$, $\Vert T_n \bigg( (t -
s)^{-1}\,\Vert x_m\Vert ^{-s}\,x_m\,{\cal
X}_{(s,t)}\bigg) \Vert < \ep$ whenever $\Vert x_m\Vert <
e^{-n/\ep}$.  Since $\Vert x_m\Vert \rightarrow 0$ as
$m \rightarrow \infty$, this is so for all $m$
sufficiently large.  It follows that
$$T_n\,\{ (t - s)^{-1}\,\Vert x_m\Vert ^{-s}\,x_m\,{\cal
X}_{(s,t)} : s < t\ ;\ m = 1, 2, 3, \ldots \}$$
is totally bounded as required.  Therefore $\{
T_n\}^\infty_{n=1}$ is a sequence of compact operators on
$Z$ such that $\Vert T_n\,z - z\Vert \rightarrow 0$ as $n
\rightarrow \infty$ for every $z$ in
$Z$.{$\spadesuit$}

\bigskip
We have in fact shown by the above argument that the
identity operator on $Z$ can be approximated, in the
topology of uniform convergence on compact sets, by
compact operators of norm at most one, that is, that $Z$
has the metric compact approximation property.

The existence of a space with the compact approximation
property but not the approximation property raises
questions as to whether results relating various stronger
versions of the approximation property, see section 1.e
of [L \& T1], have analogues for the compact
approximation property.  Some of these questions can be
answered if there is a reflexive space with the compact
approximation property but not the approximation
property, [C].  For instance, see [G \& W], Example 4.3,
where the reflexive example constructed below is used to
answer a question which arises in that paper and also in
[S].  Now $Z$ is not reflexive because the set $x_1
L^1[0,\,1]$ is a closed subspace of $Z$ which is
isomorphic to $L^1[0,\,1]$.  However, the construction of
$Z$ may be modified to produce a reflexive space as
follows.

As before, let $X$ be a Banach space which does not have
the approximation property but now suppose that $X$ is a
closed subspace of $\ell^p$ for some $2 < p < \infty$. 
This value of $p$ will remain fixed throughout the
construction.  Set $q = {p\over{p-1}}$, so that $q$ is
the conjugate of $p$.  See [L \& T1] Theorem 2.d.6 for a
proof that $\ell^p$ has such subspaces.  

Let $\{ x_n\}^\infty_{n=1}$ be a sequence in $X$ such
that: $\Vert x_n\Vert \leq 1$ for all $n$; $\Vert
x_n\Vert$ decreases to zero; and the identity operator
cannot be approximated on ${\overline{\rm
conv}}\{x_n\}^\infty_{n=1}$ by finite rank operators. 
Choose integers $n_1 < n_2 < n_3 < \ldots$ such that
$\Vert x_n\Vert < (\h)^k$ whenever $n > n_k$ and define,
for $k = 1, 2, 3,\ldots,$
$$X_k = {\rm span} \{ x_n\ :\ n \leq
n_k\}\ \ \ .$$  
Next define, for each $t$ between $0$ and $1$, $$V_t =
\bigg\{ \sum^\infty_{k=1}\,\alpha_k\,a_k/\Vert a_k\Vert ^t
\ :\ a_k \in X_k\ ,\ \Vert a_k\Vert \leq ({\scriptstyle
1\over{\raise.4ex\hbox {$\scriptstyle
2$}}})^{k-1}\ ;\
\sum^\infty_{k=1}\,\vert \alpha_k\vert ^p \leq 1\bigg\}\ \
\ .$$ 

The properties of the sets $V_t$ which we will need are
given in the following

\bigskip
\proclaim
Lemma 1.  Let $X$ be a Banach space, $X_1, X_2, X_3,
\ldots$ be finite dimensional subspaces of $X$, and $r_1,
r_2, r_3,\ldots$ be positive numbers such that
$\sum^\infty_{n=1}\,r^q_n < \infty$.  Define
$$V = \bigg\{ \sum^\infty_{k=1}\,\alpha_k\,x_k\ :\ x_k
\in X_k\ ,\ \Vert x_k\Vert \leq r_k\ ,\
\sum^\infty_{k=1}\,\vert \alpha_k\vert^p \leq 1\bigg\}\ \
\ .$$
Then $V$ is a compact, convex, symmetric subset of
$X$.  

{\sl Let $W$ be the linear subspace of $X$ spanned by $V$
and define $\Vv w \Vv = \break \inf \{ \lambda > 0\
:\ \lambda^{-1}\,w \in V\}$, $(w \in W)$.  Then
$(W,\,\Vv\cdot\Vv)$ is a Banach space. The map $Q :
(\bigoplus^\infty_{k=1}\,X_k)_{\ell^p} \rightarrow W$
defined by
$$Q(x_1,\,x_2,\,\ldots) = \sum^\infty_{k=1}\,r_k\,x_k$$
is a quotient map.}

\bigskip
\pf  It is clear that $V$ is symmetric and it is compact
because, as may be shown by a diagonal argument, every
sequence in $V$ has a convergent subsequence.  Let $x =
\sum^\infty_{k=1}\,\alpha_k\,x_k$ and $y =
\sum^\infty_{k=1}\,\beta_k\,y_k$ be in $V$.  Then
$$\equalign{\h(x + y) &= \sum^\infty_{k=1}\,\h
(\alpha_k\,x_k + \beta_k\,y_k)\cr
{}&= \sum^\infty_{k=1}\,\gamma_k\,z_k\ \ \ ,\cr}$$
where $\gamma_k = \h (\vert \alpha_k\vert + \vert
\beta_k\vert)$, $z_k$ belongs to $X_k$ and $\Vert
z_k\Vert \leq r_k$.  Since $\sum^\infty_{k=1}\,\vert
\gamma_k\vert^p \leq \sum^\infty_{k=1}\,\h (\vert
\alpha_k\vert ^p + \vert \beta_k\vert^p) \leq 1$, it
follows that $V$ is convex.  These properties of $V$
imply that $(W,\,\Vv\cdot\Vv)$ is a Banach space.

It is clear from the definition of $V$ that $Q$ maps the
unit ball of $(\bigoplus^\infty_{k=1}\,X_k)_{\ell^p}$ onto
$V$, which is the unit ball of $W$.  Therefore $Q$ is a
quotient map.{$\spadesuit$}

\bigskip
Let $W_t$ denote the space spanned by $V_t$ and normed so
that $V_t$ is its unit ball.  Denote the norm on $W_t$ by
$\Vv \cdot \Vv_t$. Then $W_s \subset W_t$ if $s < t$.  
Denote by $Q_t$ the quotient map from
$(\bigoplus^\infty_{k=1}\,X_k)_{\ell^p}$ to $W_t$ defined
by 
$$Q_t\,{\bf x} = \sum^\infty_{k=1}\,({\scriptstyle
1\over{\raise.4ex\hbox {$\scriptstyle
2$}}})^{(1-t)(k-1)}\,x_k\ \ \ \ ,$$ where ${\bf x} =
(x_1,\,x_2,\,\ldots)$.  We will require the following.

\bigskip
\proclaim
Lemma 2.  For each $s$ in $(0,\,1)$ and $\bf x$ in
$(\bigoplus^\infty_{k=1}\,X_k)_{\ell^p}$,
$$\lim_{t\rightarrow s^+}\,\Vv Q_t\,{\bf x} - Q_s\,{\bf
x}\Vv_t = 0\ \ \ .$$

\bigskip
\pf  It is clear that the limit is zero if $\bf x$ is
supported in only finitely many of the $X_k$'s and that
the set of finitely supported vectors is dense in
$(\bigoplus ^\infty_{k=1}\,X_k)_{\ell^p}$.  The result
follows because $\Vert Q_t\Vert = 1$ for each
$t$.{$\spadesuit$}

\bigskip
Let ${\cal Z}^\sharp$ be the linear span of $\{ w\,{\cal
X}_{(s,t)} : 0 < s < t < 1; w \in W_s\}$.  As in the
previous construction, ${\cal Z}^\sharp$ is a space of
functions, $f : (0,\,1) \rightarrow X$, such that $f(t)$
belongs to $W_t$ for each $t$.  Define a norm on ${\cal
Z}^\sharp$ by
$$\Vert f \Vert = \big( \int^1_0\,\Vv
f(r)\Vv^p_r\,dr\big)^{1\over p}\ \ \ ,\ \ \ (f \in {\cal
Z}^\sharp)\ \ \ .$$
Let $Z^\sharp$ denote the completion of ${\cal
Z}^\sharp$ with respect to this norm.

We will show that $Z^\sharp$ has the required
properties.  In the proofs of these properties, $\cal A$
will denote the set of all functions of the
form $$\sum^\ell_{i=1}\,\lambda_i(t_i -
s_i)^{-1/p}\,w_i\,{\cal X}_{(s_i,t_i)}\ \ \ ,\leqno(6)$$
where: $s_1 < t_1 < s_2 < t_2 < \ldots < s_\ell < t_\ell$;
$w_i \in W_{s_i}$, $\Vv w_i\Vv_{s_i} \leq 1$; and
$\sum^\ell_{i=1}\,\vert \lambda_i\vert^p = 1$.  Clearly
$\cal A$ is dense in the unit ball of $Z^\sharp$.

\bigskip
\proclaim
Proposition 3.  $Z^\sharp$ is a quotient of a closed
subspace of $L^p(0,\,1)$.  In particular, $Z^\sharp$ is
reflexive.

\bigskip
\pf  Let $X_k$, $k = 1,\,2,\,3,\,\ldots$ be the subspaces
of $\ell^p$ defined above.  Then
$(\bigoplus^\infty_{k=1}\,X_k)_{\ell^p}$ is a subspace of
$(\bigoplus^\infty_{k=1}\,\ell^p)_{\ell^p}$.  It follows
that $L^p\big((0,\,1)\,,\,(\bigoplus^\infty_{k=1}
\,X_k)_{\ell^p}\big)$ is isometric to a subspace of
$L^p(0,\,1)$.  We will show that $Z^\sharp$ is a quotient
of this space.

The quotient map $Q :
L^p\big((0,\,1)\,,\,(\bigoplus^\infty_{k=1}\,X_k)_{\ell^p}
\big) \rightarrow Z^\sharp$ will be defined first of all 
on simple functions from $(0,\,1)$ to
$(\bigoplus^\infty_{k=1}\,X_k)_{\ell^p}$.  Let $f$ be such
a function and put
$$(Q\,f)(t) = Q_t(f\,(t))\qquad,\qquad (0 < t < 1)\ \ \
.\leqno(7)$$
Then $(Q\,f)(t)$ belongs to $W_t$ for each $t$ and it
follows from Lemma 2 that there are functions $f_n$, $n =
1,\,2,\,3,\,\ldots$ in ${\rm span}\{ w\,{\cal X}_{(s,t)}
: 0 < s < t < 1\ \ ;\ \ w \in W_s\}$ such that
$$\lim_{n\rightarrow\infty}\,( \int^1_0\,\Vv (Q\,f)(r) -
f_n(r) \Vv^p_r\,dr)^{1\over p} = 0\ \ \ .$$
It follows that: $\{f_n\}^\infty_{n=1}$ is a Cauchy
sequence in $Z^\sharp$;  $Q\,f$ may be identified with
the limit of this sequence;  and
$$\equalign{\Vert Q\,f\Vert &= (\int^1_0\,\Vv
(Q\,f)(r)\Vv^p_r\,dr)^{1\over p}\cr
&= (\int ^1_0\,\Vv Q_r(f(r))\Vv^p_r\,dr)^{1\over p}\ \ \
.\cr}$$
Since $\Vert Q_r\Vert \leq 1$ for each $r$, $\Vert Q\Vert
\leq 1$.  Since the simple functions are dense, $Q$
extends by continuity to be a norm one operator from
$L^p\big((0,\,1)\ \ ,\ \
(\bigoplus^\infty_{k=1}\,X_k)_{\ell^p}\big)$ to
$Z^\sharp$.  The equation (7) will hold for every $f$ in
$L^p\big((0,\,1)\ \ ,\ \
(\bigoplus^\infty_{k=1}\,X_k)_{\ell^p}\big)$ and almost
every $t$ in $(0,\,1)$.

Now let $f = \sum^\ell_{i=1}\,\lambda_i(t_i -
s_i)^{-1/p}\,w_i\,{\cal X}_{(s_i,t_i)}$ be in $\cal A$,
where $w_i = \sum^\infty_{k=1}\,\alpha_k\,a_k/\Vert
a_k\Vert^{s_i}$ belongs to $V_{s_i}$ for each $i$.  For
each $s_i \leq t \leq t_i$, put 
$${\bf w}_i^{(t)} = (\alpha_1\,a_1/ \Vert
a_1\Vert^{s_i}, \ldots ,\alpha_k\,({\scriptstyle
1\over{\raise.4ex\hbox {$\scriptstyle 2$}}})
^{(t-1)(k-1)}a_k/ \Vert a_k\Vert^{s_i}\,,\ldots)$$ in
$(\bigoplus^\infty_{k=1}\,X_k)_{\ell^p}$.  Then $\Vert
{\bf w}^{(t)}_i\Vert \leq\big(\sum^\infty_{k=1}\,\vert
\alpha_k\vert^p\, ({\scriptstyle 1\over{\raise.4ex\hbox
{$\scriptstyle 2$}}}) ^{p(t-s_i)(k-1)}\big)^{1\over p} \leq
1$ and $Q_t ({\bf w}_i^{(t)}) = w_i$ for each $s_i \leq t
\leq t_i$.  Hence, if we define $${\bf f}(t) =
\sum^\ell_{i=1}\,\lambda_i(t_i - s_i)^{-1/p}\,{\bf
w}_i^{(t)}\,{\cal X}_{(s_i,t_i)}(t)\qquad ,\qquad (0 < t <
1)\ \ \ ,$$ then $\Vert {\bf f}\Vert \leq 1$ and $Q{\bf f}
= f$.  Therefore $Q$ maps the unit ball of
$L^p\big((0,\,1)\ \ ,\ \
(\bigoplus^\infty_{k=1}\,X_k)_{\ell^p}\big)$ onto the
unit ball of $Z^\sharp$ and so $Q$ is a quotient
map.{$\spadesuit$}

Since $L^p(0,\,1)$ is uniformly convex, $L^p\big((0,\,1)\
,\ (\bigoplus^\infty_{k=1}\,X_k)_{\ell^p}\big)$ is
uniformly convex and so it follows, by [Da] Theorem 5.5,
that $Z^\sharp$ is also uniformly convex.  This theorem
of Day may also be deduced from the duality between uniform
convexity and uniform smoothness, (see [L \& T2],
Proposition 1.e.2).

\bigskip
\proclaim
Proposition 4.  $Z^\sharp$ has the metric compact
approximation property but does not have the
approximation property.

\bigskip
\pf  The same argument as showed that $Z$ does not have
the approximation property also shows that $Z^\sharp$
does not have this property.  It is immediate from the
definitions of these sets that $U_{1\over 2} \subseteq
V_{1\over 2}$ and so $Y_{1\over 2}$ is contained in
$W_{1\over 2}$ and the embedding is a contraction. 
Hence maps $R^\sharp : Y_{1\over 2} \rightarrow Z^\sharp$
and $J^\sharp : Z^\sharp \rightarrow X$ are defined by
equations (1) and (3) respectively and an estimate
similar to (2) shows that $R^\sharp$ is compact.  These
new maps also satisfy that $J^\sharp\,R^\sharp = L_{1\over
2}$ and so $Z^\sharp$ does not have the approximation
property.

Shift operators $S^\sharp_r : Z^\sharp \rightarrow
Z^\sharp$ may be defined for each $r$ between $0$ and $1$
just as they were on $Z$ and then used to define
operators $T^\sharp_n : Z^\sharp \rightarrow Z^\sharp$
for $n = 1, 2, 3, \ldots$ by an equation similar to (4). 
The same argument as used in Proposition 2 shows that
these operators have norm at most one and that $\{
T^\sharp_n\}^\infty_{n=1}$ converges to the identity
operator uniformly on compact subsets of $Z^\sharp$. 
However, to show that $T^\sharp_n$ is compact for each
$n$ requires a slightly different argument to that used
in Proposition 2. It suffices to show that
$T^\sharp_n\,{\cal A}$ is totally bounded.

For each positive integer $m$ define the subset,
$V^{(m)}_s$, of $V_s$ by
$$V_s^{(m)} = \bigg\{
\sum^\infty_{k=m}\,\alpha_k\,a_k/\Vert a_k\Vert^s\ :\ a_k
\in X_k\ ,\ \Vert a_k\Vert \leq ({\scriptstyle
1\over{\raise.4ex\hbox {$\scriptstyle
2$}}})^{k-1}\ ;\ \sum^\infty_{k=m}\,\vert
\alpha_k\vert ^p \leq 1\bigg\}\ \ \ .$$
Let $w = \sum^\infty_{k=m}\,\alpha_k\,a_k/\Vert
a_k\Vert^s$ be in $V^{(m)}_s$ and suppose that $\Vv
w\Vv_s = 1$.  Then, for $r > s$,
$$\equalign{\Vv w\Vv_r &= \Vv
\sum^\infty_{k=m}\,\alpha_k\,\Vert
a_k\Vert^{r-s}\,a_k/\Vert a_k\Vert^r\,\Vv _r\cr
{}&\leq \bigg( \sum^\infty_{k=m}\,\vert \alpha _k\,\Vert
a_k\Vert^{r-s}\,\vert^p\bigg)^{1\over p}\cr
{}&\leq ({\scriptstyle 1\over{\raise.4ex\hbox
{$\scriptstyle 2$}}})^{(m-1)(r-s)}\ \ \ .\cr}$$
Corresponding to the estimate in (5) we now have
$$\lequalignno{\Vert T^\sharp_n (w\,{\cal X}_{(s,t)})\Vert
&\leq n (t - s)\,\left(\int^{t+1/n}_s\,\Vv
w\Vv^p_r\,dr\right) ^{1\over p}&(8)\cr
{}&< n(t - s)(m - 1)^{-{1\over p}}\ \ \
.\cr}$$ Next, for each positive integer $m$
denote by ${\cal B}^{(m)}$ the subset of $\cal
A$ consisting of all functions of the form (6)
where $w_i$ belongs to $V^{(m)}_s$ for each
$i$ and by ${\cal C}^{(m)}$ the subset of
$\cal A$ consisting of functions where $w_i$
belongs to the finite dimensional space ${\rm
span}\{x_n : n \leq n_{m-1}\} = X_{m-1}$.  For a general
$w$ in $V_s$, $w = \sum^{m-1}_{k=1}\, \alpha_k\,a_k/\Vert
a_k\Vert^s + w^{\prime}$, where $w^{\prime}$ is in
$V^{(m)}_s$ and $\sum^{m-1}_{k=1}\,\alpha_k\,a_k/\Vert
a_k\Vert^s$ belongs to the unit ball of $X_{m-1}$.  It
follows that ${\cal A} \subseteq {\cal B}^{(m)} +
{\cal C}^{(m)}$.

For each $\xi$ in ${\cal B}^{(m)}$ we
have, by (6) and (8),
$$\equalign{\Vert T^\sharp_n\,\xi\Vert
&\leq n(m - 1)^{-{1\over p}}\,\sum^\ell_{i=1}\,\vert
\lambda_i\vert\ \vert t_i - s_i\vert^{1-{1\over p}}\cr
{}&\leq n(m - 1)^{-{1\over p}}\ \ \ ,\cr}$$
because $\sum^\ell_{i=1}\,\vert\lambda_i\vert^p = 1$ and
$\sum ^\ell_{i=1}\,\vert t_i - s_i\vert \leq 1$. If, given
$\ep > 0$, we choose $m > 1 + (n/\ep)^p$, then it follows
that $\Vert T^\sharp_n\,f\Vert < \ep$ for every $f$ in
${\cal B}^{(m)}$.  Also, since the map $f \mapsto nf *
{\cal X}_{(0,1/n)}$ is a compact operator on $L^p(0,\,1)$,
$T^\sharp_n\,{\cal C}^{(m)}$ is totally bounded for each
$m$.  Therefore $T^\sharp_n{\cal A}$ is totally bounded
and so $T^\sharp_n$ is a compact operator.$\spadesuit$

It may also be shown that, for $1 < p < 2$, there are
quotients of subspaces of $L^p(0,\,1)$ which have the
metric compact approximation property but not the
approximation property.  This may be shown by choosing a
subspace, $X$, of $\ell^p$ which does not have the
approximation property, such spaces have been shown to
exist by Szankowski, see [Sz] or [L\&T2], theorem 1.g.4,
and then repeating the above construction. 
Alternatively, the dual of the above example has the
required properties.  I am grateful to Professor T. Figiel
for this remark and also for some other suggestions which
shortened some proofs and improved Proposition 3.

\vskip2.5in
Part of this work was completed while the author was
visiting the University of Leeds with the generous
support of SERC grant GR-F-74332.  It is a pleasure to
thank the members of the School of Mathematics at Leeds
for their hospitality and to thank Michel Solovej and
Niels Gr{\o}nb{\ae}k for helpful conversations concerning
the subject of this paper. This paper is accepted for 
publication in Studia Mathematica.

\vfill\eject
\centreline{\bf References}

\bigskip\bigskip
\ref
[C]\quad Casazza, P., personal communication

\bigskip
\ref
[D1]\quad Davie, A.M., ``The approximation problem for
Banach spaces'', Bull. London Math. Soc. {\bf 5}, (1973),
261-266.

\bigskip
\ref
[D2]\quad Davie, A.M., ``The Banach approximation
problem'', J. Approx. Theory {\bf 13}, (1975), 392-394.

\bigskip
\ref
[Da]\quad	Day, M.M., ``Uniform convexity in factor and
conjugate spaces'', Ann. of Maths. {\bf 45}, (1944),
375-385.

\bigskip
\ref
[E]\quad Enflo, P., ``A counterexample to the
approximation property in Banach spaces'', Acta Math.
{\bf 130}, (1973), 309-317.

\bigskip
\ref
[G]\quad Grothendieck, A., ``Produits tensoriels
topologiques et espaces nucleaires'', Mem. Amer. Math.
Soc. {\bf 16} (1955).

\bigskip
\ref
[G \& W]\quad Gr{\o}nb{\ae}k, N. and Willis, G., 
``Approximate identities in Banach algebras of compact
operators'', to appear in Can. Math. Bull.

\bigskip
\ref
[L \& T1]\quad Lindenstrauss, J. and Tzafriri, L.,
Classical Banach Spaces I, Berlin-Heidelberg-New York,
Springer 1977.

\bigskip
\ref
[L \& T2]\quad Lindenstrauss, J. and Tzafriri, L.,
Classical Banach Spaces II, Berlin-Heidelberg-New York,
Springer 1979.

\bigskip
\ref
[S]\quad Samuel, C., ``Bounded approximate identities in
the algebra of compact operators on a Banach space'', to
appear in Proc. Amer. Math. Soc.

\end
\def\ep{\varepsilon} 
\def\ch{\raise.4ex\hbox{$\chi$}}
\def\va{\raise.4ex\hbox{$\varphi$}}
\def\centreline{\centerline} \def\equalign{\eqalign}
  \def\lequalignno{\leqalignno}
\def\nt{\noindent} 
 \def\ga{\raise.4ex\hbox{$\gamma$}}

\def\Vv{{|||}}   \def\pf{\noindent{\bf Proof.\ \ }} 
    
\font\titlefont=cmbx10 scaled\magstep 3
\font\namefont=cmbx10 scaled\magstep 1
\def\ref{\noindent\hang}
\def\SA{\lower.9ex\hbox{$\buildrel \subset \over \rightarrow$}}
\def\BTD{\raise.4ex\hbox{$\bigtriangledown$}}

\def\MATHOP(#1){\mathop{#1}}

\def\AST{\lower.8ex\hbox{*}}

\def\h{{1\over 2}}
 3
 1
\def\Rogerslash{\backslash\kern-2pt\backslash}
\def\RS{/\kern-2pt /}

\def\GH{G{\mathop{\times
\atop{\raise.4ex\hbox{\sevenrm{\ \ s.d.} }}}}H}




\def\vz{\ |\kern-3pt-\kern-5pt-\ }
\def\vr{\ |\kern-2pt-\kern-5pt\rightarrow\ }

\def\ontoarrow{\geq\kern-3pt-\kern-5pt-\kern\gg}


\def\bookformat{ 
\hsize=157true mm\vsize=240true mm\parindent=25true mm
\baselineskip=18true pt \voffset=36true pt \nopagenumbers
\def\rhead{\hfil} \def\lhead{\hfil}
 \headline={\sevenrm\it\ifodd\pageno\rhead\else\lhead\fi}
\def\newsection##1{\vskip20true pt{\bf##1}
\def\rhead{\centreline{##1}\rlap{\mit\folio}}
\def\lhead{\llap{\mit\folio}\centreline{##1}}}
\def\newchapter##1{{\bf##1}\vskip53true mm}}

\def\newleftmargin#1mm#2mm{
\def\funnymargin{\ifodd\pageno#1mm\else#2mm\fi}
\hoffset=\funnymargin}
 
\def\onto{\ >\kern-5pt-\kern-2pt\kern-4pt\gg\ }